\documentclass[12pt]{article}

\usepackage{amsmath}
\usepackage{amssymb}
\usepackage{bbold}
\usepackage{color}
\usepackage{graphicx}
\usepackage{amsthm}
\usepackage{multirow}
\usepackage{threeparttable}

\newtheorem{theorem}{Theorem}
\newtheorem{corollary}{Corollary}
\theoremstyle{definition}

\theoremstyle{remark}
\newtheorem{remark}{Remark}

\title{Quantum hub and authority centrality measures for directed networks based on continuous-time quantum walks}

\author{P.~Boito\footnote{Universit\`a di Pisa, Dipartimento di Matematica, Largo Bruno Pontecorvo 5, 56127 Pisa, Italy. E-mail: {\tt paola.boito@unipi.it}} \mbox{}  and R.~Grena\footnote{C. R. ENEA Casaccia, via Anguillarese 301, 00123 Roma, Italy. E-mail: {\tt roberto.grena@enea.it}}}
\date{}

\begin{document}

\maketitle

\section*{Abstract}

In this work we introduce, test and discuss three quantum methods for computing hub and authority centrality scores in directed networks. The methods are based on unitary, continuous-time quantum walks; the construction of a suitable Hermitian Hamiltonian is achieved by performing a quantum walk on the associated bipartite graph. 

Two methods, called CQAu and CQAw, use the same evolution operator, inspired by the classical HITS algorithm, but with different initial states; the computation of hub and authority scores is performed simultaneously. The third method, called CQG and inspired by classical PageRank, requires instead two separate runs with different evolution operators, one for hub and one for authority scores. 

The methods are tested on several directed graphs with different sizes and properties; a comparison with other well-established ranking algorithms is provided. CQAw emerges as the most reliable of the three methods and yields rankings that are largely compatible with results from HITS, although CQAu and CQG also present interesting features and potential for applications.

\section{Introduction}\label{sec:intro}

The design and analysis of centrality measures for networks -- i.e., measures of importance of nodes -- plays a crucial role in network analysis and has sparked much interest in the last decades; a classical example is the Page Rank algorithm \cite{BrinPage} for ranking web pages in search engines. 

For directed networks, moreover, there is a distinction to be made between {\em hub} and {\em authority} centrality. For instance, in a network of Web pages or documents, hubs can be seen as pointers to high quality resources, whereas authority are the resources themselves. 
Therefore each node plays a double role in the network and is assigned separate hub and authority scores. This notion was initially popularized by the HITS algorithm \cite{hits} for ranking Web pages (see Section \ref{sec:background}). Benzi et al. \cite{BEK} have proposed a bipartization approach where hub and authority centrality is defined through the exponential of the adjacency matrix of the undirected bipartite graph associated with the original graph.  

Meanwhile, interest in quantum computation fueled the development of the theory of quantum walks on networks, the quantum version of classical random walks. Quantum walks have, of course, the advantage of exploiting the considerable potential of quantum algorithms for fast computation, and exhibit a few peculiar differences w.r.t. classical quantum walks. A good introduction to quantum walks can be found, for instance, in the recent book \cite{Portugal}; we mention also the seminal paper \cite{Aharonov} and the review \cite{Venegas}.

As in the case of classical random walks, quantum walks are available in two flavors: discrete-time quantum walks (DTQW) and continuous-time quantum walks (CTQW). DTQWs typically take place in a Hilbert space with dimension higher than the number of nodes, because extra degrees of freedom are added 
in order to obtain a unitary evolution operator $U$ modeling each step. In the coin-flip formalism, for instance, the dimension of the Hilbert space is of the order of the number of edges in the graph. In CTQW, on the other hand, the underlying Hilbert space usually has dimension equal to the number of nodes and the evolution operator is $U(t) = \exp(-itH)$, where $H$ is a suitably defined Hermitian Hamiltonian operator that encodes the structure of the network.

Following the idea behind PageRank, where centrality is defined as the stationary vector of a (classical) random walk, several authors have proposed quantum algorithms for centrality computation based on quantum walks. A DTQW-based ranking algorithm, called {\em Quantum PageRank} (DQPR), was proposed in \cite{Paparo}. The dimension of the Hilbert space used in DQPR is $n^2$, where $n$ is the number of nodes. Another DTQW-based approach is presented in \cite{Chawla}. The use of DTQW has also been proposed in connection with community detection in networks \cite{MukaiHatano}. Unitary CTQW have been applied for ranking purposes in undirected networks \cite{Rossi, Izaac}: in this case the quantum walk takes place in a Hilbert space of dimension $n$. An extension of the idea in \cite{Izaac} to directed networks forgoes the requirement for a unitary evolution operator and adopts instead non-standard PT-symmetric Hamiltonians \cite{IzaacPT}. See also \cite{Wang} for a physical implementation of a CTQW-based approach to ranking in directed networks.

A common feature of unitary quantum walks, due to the unitary evolution, is the absence of asymptotic states. As a consequence, 
the quantum state of the system (which encodes the occupation values of each node) does not converge to a stationary state: the extraction of centrality scores usually requires step- or time-averages. Another consequence is that the occupation of each node may depend on the initial state, i.e., the prescription of the initial state is part of the ranking method. We mention however that there are also quantum methods that obtain stable states without averaging, through a tunable time-dependent Hamiltonian (adiabatic computations, as in \cite{PRadiabatic}) or using mixed quantum-classical evolution \cite{PRsemiclass}. 

In this work we introduce and discuss three new unitary-CTQW-based hub and authority centrality measures for directed networks. Two different Hamiltonians are used: one based directly on the adjacency matrix, the other one built starting from the Google matrix. We rely on bipartization in order to define Hermitian Hamiltonians, and consequently unitary evolution operators. The dimension of the underlying Hilbert space is $2n$ (that is, the number of nodes of the bipartite graph associated with the given network). In two of the three methods, hub and authority rankings are computed simultaneously.

Section \ref{sec:background} provides the background to understand the problem and briefly presents a few popular ranking algorithms that will later be used for comparison purposes. Our CTQW-based methods are outlined in Section \ref{sec:main} according to the evolution operators they use. Section \ref{sec:experiments} is devoted to numerical experiments: tests on simple toy models, which allow us to discuss certain features of the methods, tests on artificial 128-nodes scale-free graphs, and tests on real-life larger graphs (up to 4472 nodes). The results are compared to rankings obtained from four other algorithms, three classical (HITS, PageRank and the method in \cite{BEK}) and one quantum (QPR). Section \ref{sec:conclusions} discusses the proposed methods in the light of the numerical results.

\section{Background}\label{sec:background}

For our purposes, a {\em (directed) graph} ${\mathcal G}$ is a pair $(V,E)$, where $V$ is the set of nodes, labeled from $1$ to $n$, and $E=\{(i,j)|i,j\in V\}$ is the set of edges. The graph is called {\em undirected} if $(i,j)\in E$ implies $(j,i)\in E$. The graphs we consider are weakly connected, unweighted and contain no loops (i.e., edges from one node to itself) or multiple edges.

Recall that the {\em adjacency matrix} associated with a graph ${\mathcal G}$ is the $n\times n$ matrix $A=(A_{ij})$ such that $A_{ij}=1$ if there is an edge from node $i$ to node $j$, $A_{ij}=0$ otherwise. Clearly $A$ is symmetric if and only if ${\mathcal G}$ is undirected.

The {\em in-degree} ${\rm deg}_{\rm in}(i)$ of node $i$ is the number of nodes pointing to node $i$, whereas its {\em out-degree} ${\rm deg}_{\rm out}(i)$ is the number of nodes pointed to by node $i$. In matrix terms, it holds ${\rm deg}_{\rm in}(i)=({\mathbb{1}}^T A)_i$ and ${\rm deg}_{\rm out}(i)=(A {\mathbb{1}})_i$, where ${\mathbb{1}}$ is the column vector formed by $n$ components all equal to $1$.

A graph is called {\em bipartite} if its vertices can be divided into two disjoint sets, such that no edge of the graph connects vertices from the same set. Given a directed graph ${\mathcal G}=(V,E)$, an associated undirected bipartite graph ${\mathcal B}_{\mathcal G}$ can be defined as follows:
\begin{itemize} 
\item take as set of nodes the union of two copies $V'$ and $V''$ of $V$, 
\item nodes $i'\in V'$ and $j''\in V''$, corresponding to nodes $i$ and $j$ in $\mathcal{G}$, are connected by an edge if and only if there is an edge from $i$ to $j$ in ${\mathcal G}$.
As the graph is bipartite, there is no edge between nodes of $V'$, nor between nodes of $V''$.
\end{itemize}  
See Figure \ref{fig:BEK} for an example of bipartite graph associated with a given directed graph.  

Note that, if the nodes of ${\mathcal B}_{\mathcal G}$ are labeled by listing first the nodes in $V'$ and then the nodes in $V''$, in the same order as in $V$, then the adjacency matrix of ${\mathcal B}_{\mathcal G}$ takes the form 
$\left[\begin{array}{cc} 0 & A\\ A^T & 0\end{array}\right],$ where $A$ is the adjacency matrix of ${\mathcal G}$.

\subsection{Classical centrality measures for directed graphs}

Among the numerous centrality measures proposed in the literature for directed graphs, we recall a few that are especially popular and that will be used for comparison purposes in Section \ref{sec:experiments}. Two of them (HITS and PageRank) were originally developed for ranking web pages.

The HITS (Hyperlink-Induced Topic Search) algorithm \cite{hits} is an iterative scheme that computes an authority score and a hub score for each node. Each iteration takes the form
\begin{equation}
{ x}^{(k)}=A^T { y}^{(k-1)},\qquad { y}^{(k)}=A { x}^{(k)},\label{eq:hits}
\end{equation}
followed by normalization in $2$-norm, where ${ x}^{(k)}$ and ${ y}^{(k)}$ are the vectors of authority and hub scores, respectively, at the $k$-th step. Authority and hub centralities are defined as the limit values of the corresponding scores for $k\rightarrow\infty$. This process can also be reframed as a power method yielding the dominant eigenvector of $A^T A$ (for authorities) and of $A A^T$ (for hubs).
As the Perron-Frobenius theorem does not necessarily apply to these matrices, the scores computed by HITS may depend on the choice of the initial vector. See \cite{Farahat} for an analysis of existence and uniqueness properties of the HITS vectors.

PageRank \cite{BrinPage, LangvilleMeyer} is very well known, following its success as a web ranking algorithm. In the PageRank algorithm, the adjacency matrix $A$ of the graph is modified to obtain the so-called patched adjacency matrix $\tilde{A}$, which is row-stochastic:
$$
\tilde{A}_{ij}=\left\{\begin{array}{ll}
A_{ij}/{\rm deg_{\rm out}}(i) & {\rm if }\, {\rm deg_{\rm out}}(i)>0,\\
1/n & {\rm otherwise.}
\end{array}\right.
$$
Each entry $\tilde{A}_{ij}$ can be seen as the probability, for a walker placed at node $i$, to reach node $j$, with the assumption that from each dangling node -- i.e., a node with zero out-degree -- the walker can reach all other nodes with uniform probability (or, more generally, with a probability given by a ``personalization vector''). Next, one introduces a ``teleportation'' or ``random surfer'' effect by adding a rank-one correction:
\begin{equation}
G=\alpha\tilde{A}+\frac{1-\alpha}{n} {\mathbb 1} {\mathbb 1}^T,\label{eq:googlematrix}
\end{equation} 
where $\alpha\in[0,1]$ is a parameter, usually chosen as $\alpha=0.85$. The matrix $G$ is sometimes referred to as the Google matrix. By the Perron-Frobenius theorem, $1$ is the (simple) eigenvalue of maximum modulus for $G^T$. The associated eigenvector $x>0$ such that
$G^Tx=x$ yields the PageRank scores of the nodes, and it can be computed efficiently via the power method.

Note that PageRank defines authority scores, since the importance of each node is determined by the importance of the nodes pointing to it. Conversely, {\em reverse PageRank} \cite{Fogaras}, that is, PageRank applied to the ``reversed graph'' with adjacency matrix $A^T$, provides hub scores.

The third approach we mention is presented in \cite{BEK} and can be placed in the more general framework of walk-based centrality mesures defined through matrix functions. Recall that the notion of {\em subgraph centrality} for undirected graphs takes as centrality scores the diagonal entries of ${\rm e}^A$, where $A$ is the (symmetric) adjacency matrix \cite{EstradaHigham}. In the directed case, the authors of \cite{BEK} consider instead the exponential of the symmetric $2n\times 2n$ augmented matrix
\begin{equation}
\mathcal{A}=\left[\begin{array}{cc} 0 & A \\ A^T & 0\end{array}\right].\label{augmented}
\end{equation}
Observe that, as mentioned above, $\mathcal{A}$ is the adjacency matrix of the undirected bipartite graph associated with the directed graph under consideration. The hub and authority scores of the $i$-th node are defined as $[{\rm e}^{\mathcal{A}}]_{i,i}$ and $[{\rm e}^{\mathcal{A}}]_{n+i,n+i}$, respectively. 
A walk-based interpretation of these quantities relies on the observation that $[(AA^T)^k]_{ij}$ counts the number of alternating walks of length $2k$ that connect node $i$ to node $j$ starting with an out-edge, whereas $[(A^TA)^k]_{ij}$ is the number of alternating walks of length $2k$ from node $i$ to node $j$ starting with an in-edge. 

\subsection{Quantum centrality}

Quantum walks (i.e., quantum versions of classical random walks) have been used by several authors to define measures of centrality on graphs.
A detailed presentation of quantum walks goes beyond the scope of the present work; we refer the interested reader to the seminal paper \cite{Aharonov} and to the review \cite{Venegas}, as well as the book \cite{Portugal}.

Much of the literature on quantum centrality relies on discrete-time quantum walks: one notable example is the (discrete) Quantum Page Rank, here abbreviated as DQPR \cite{Paparo}. Note that, for a graph with $n$ nodes, these methods typically work in an $n^2$-dimensional state space, as required by the coined or Szegedy formulations of the DTQW (although in some cases, as for DQPR, it can be proved that the non-trivial component of the dynamics actually takes place in a subspace of dimension at most $2n$).

Continuous-time quantum walks \cite{Feynman, FarhiGutmann, Childs} have also been used to define notions of quantum centrality for undirected networks. 
A CTQW on a graph ${\mathcal G}$ with $n$ vertices is a process that takes place in a complex Hilbert space of dimension $n$, spanned by basis states $|j\rangle$, with $j=1,\ldots,n$, each of them associated with the corresponding node. A generic state of the system is denoted as $|\psi\rangle$ and takes the form
$$
|\psi\rangle=\sum_{j=1}^n a_j |j\rangle,
$$ 
where the amplitudes $a_1,\ldots,a_n\in\mathbb{C}$ are such that $\sum_{j=1}^n |a_j|^2=1$. Note that each amplitude $a_j$ can be expressed as the scalar product $\langle j | \psi\rangle$. Its square modulus $|a_j|^2$ is the {\em occupation} of node $j$, i.e., the probability of finding the system in the state $|j\rangle$ after a measurement.

The time evolution of a CTQW on a graph ${\mathcal G}$ is described by the Schr\"odinger equation
\begin{equation}
i\frac{\partial |\psi(t)\rangle}{\partial t} = H |\psi(t)\rangle, \label{eq:schroedinger}
\end{equation}
where $|\psi(t)\rangle$ is the state of the system at time $t$, and the Hamiltonian operator $H$ encodes the structure of the graph. Two common choices for the Hamiltonian are the graph Laplacian and the adjacency matrix, which are symmetric in the hypothesis of undirected ${\mathcal G}$; see \cite{Wong} for motivation and a comparison. 

The general solution of \eqref{eq:schroedinger} yields the state of the system at time $t$ and is given by 
\begin{equation}
|\psi(t)\rangle=U(t) |\psi(0)\rangle,\label{eq:solution}
\end{equation}
where $U(t)={\rm e}^{-iHt}$ is the unitary evolution operator of the system and $|\psi(0)\rangle$ is the initial state. Quantum properties of the CTQW include superposition of states and time reversibility; the latter is a consequence of the unitarity of the evolution operator and implies that there is no limiting state. 

We recall here a CTQW-based centrality measure introduced in \cite{Izaac} for undirected graphs. This approach has the merit of limiting the dimension of the Hilbert space to $n$, for an $n$-node graph. It can be summarised in three main steps:
\begin{enumerate}
\item Prepare the initial state (the ``walker'') in an equal superposition of all vertex states:
$$|\psi(0)\rangle=\frac{1}{\sqrt{n}}\sum_{j=1}^n|j\rangle.$$
\item Propagate the walker for a sufficiently large time $t$:
$$|\psi(t)\rangle={\rm e}^{-iHt}|\psi(0)\rangle,$$ where the Hamiltonian $H$ is chosen as the adjacency matrix of the graph.
\item Compute the time-average probability distribution of finding the walker at each vertex $j$; this is the centrality score of node $j$:
$$C_j=\lim_{t\rightarrow\infty}\frac{1}{t}\int_0^t|\langle j|\psi(t)\rangle|^2 dt.$$
\end{enumerate}
Note that the introduction of time averages at step 3 is made necessary by the absence of a limit state and is a recurrent ingredient in the analysis of both discrete- and continuous-time quantum walks.

An interesting attempt at modeling quantum centrality on directed graphs through a nonunitary evolution operator has been proposed in \cite{IzaacPT} and \cite{Wu} exploiting parity-time symmetry.

\section{Three centrality methods in directed networks based on unitary CTQW}\label{sec:main}

Building on ideas from \cite{Izaac} and \cite{BEK}, in this work we propose to employ continuous-time quantum walks to define a measure of hub and authority centrality for directed graphs, while maintaining unitarity.
We apply a bipartization approach: broadly speaking, one can think of obtaining hub and authority centrality scores for ${\mathcal G}$ by applying an undirected centrality measure to the bipartite graph associated with ${\mathcal G}$. This idea is put forward in \cite{BEK} for the exponential-based subgraph centrality. It also holds for HITS: one may equivalently apply HITS to a directed graph and compute the hub and authority scores $x$ and $y$ as in \eqref{eq:hits}, or consider the associated undirected bipartite graph, compute the positive dominant eigenvector $u$ of ${\mathcal A}$ and take $u(1:n)$ as hub scores, $u(n+1:2n)$ as authority scores (see again \cite{BEK}). In other words, we can think of nodes $1,\ldots,n$ in the bipartite graph as corresponding to the nodes of ${\mathcal G}$ in their role as hubs, and of nodes $n+1,\ldots,2n$ as the nodes of ${\mathcal G}$ in their role as authorities.

\subsection{Two methods based on the adjacency matrix: CQAu and CQAw}

Now, let ${\mathcal G}$ be a directed graph with $n$ nodes and let $A$ be the associated adjacency matrix. As in PageRank, consider the rank-one correction
\begin{equation}
\tilde{A}=\alpha A+\frac{1-\alpha}{n}{\mathbb{1}}{\mathbb{1}}^T,\label{eq:rank1}
\end{equation}
where $\alpha\in [0,1]$ is a suitably chosen parameter. Define the symmetric matrix $H$ of size $2n\times 2n$ as follows:
\begin{equation}
H=\left[\begin{array}{cc}
0 & \tilde{A} \\ \tilde{A}^T & 0
\end{array}\right].\label{eq:hamiltonian}
\end{equation}
The matrix $H$ represents our Hamiltonian operator, which will be used to model a CTQW in a Hilbert space of dimension $2n$.

The motivation for introducing the correction \eqref{eq:rank1} stems from the following considerations. If $\alpha=1$, that is, if there is no rank-one correction of the adjacency matrix, then the Hamiltonian $H$ is the adjacency matrix of the undirected bipartite graph ${\mathcal B}_{\mathcal G}$ associated with ${\mathcal G}$. 
Note that ${\mathcal B}_{\mathcal G}$ might not be connected, even if $\mathcal G$ is strongly connected; see, e.g., Example 5 in Section \ref{sec:experiments}.
 In such cases, the Hamiltonian $H$ with $\alpha=1$ gives rise to a biased model, because the resulting CTQW is ``trapped'' in separate connected components. An extreme example of such behavior arises in presence of a node with zero in-degree or out-degree: one of its corresponding nodes in ${\mathcal B}_{\mathcal G}$ turns out to be disconnected from the rest of the graph, so its occupation remains fixed at the initial value throughout the walk.
The rank-one correction in \eqref{eq:rank1} solves this difficulty and adds a moderate random surfer effect. In our experiments we take $\alpha=0.85$ as is often done for classical PageRank.
  
The CTQW is now described by the Schr\"odinger differential equation \eqref{eq:schroedinger} with $H$ as in \eqref{eq:hamiltonian}. Therefore the state of the system at time $t$ is
$$
|\psi(t)\rangle=U(t) |\psi(0)\rangle,
$$
where the evolution operator $U(t)={\rm e}^{-iHt}$ is the exponential of a Hermitian matrix and is therefore unitary. (More details on the choice of the initial state $|\psi(0)\rangle$ will be given later in this section).

With this approach, the hub centrality of node $j$ in $\mathcal{G}$ is defined as
\begin{equation}
C_{\rm hub}(j)=\lim_{T\rightarrow\infty}\frac{1}{T}\int_0^T|\langle j|\psi(t)\rangle|^2 dt \label{eq:ourhub}
\end{equation}
and the authority centrality of node $j$ is
\begin{equation}
C_{\rm auth}(j)=\lim_{T\rightarrow\infty}\frac{1}{T}\int_0^T|\langle j+n|\psi(t)\rangle|^2 dt. \label{eq:ourauthority}
\end{equation}
The limits in \eqref{eq:ourhub} and \eqref{eq:ourauthority} are always well-defined. In fact we have the following result, adapted from its discrete-time version in \cite{Aharonov}:
\begin{theorem}\label{th:Thm1}
Let $U(t)={\rm e}^{-itH}$ be the unitary evolution operator of the CTQW and denote as $\{\theta_j, |\phi_j\rangle\}_{j=1,\ldots,N}$, with $N=2n$, the eigenvalues and eigenstates of the Hamiltonian $H$. Recall that the eigenvalues and eigenstates of $U(t)$ are $\{\lambda_j(t)= {\rm e}^{-i t\theta_j}, |\phi_j\rangle\}_{j=1,\ldots,N}$.
  Let $|\psi(0)\rangle=\sum_{i=1}^N a_j|\phi_j\rangle$, with $a_1,\ldots,a_N\in\mathbb{C}$, be the initial state of the system, written in the eigenstate basis, and denote as $|\psi(t)\rangle$ the state of the system at time $t$. Then for any state $|\xi\rangle$ it holds
$$
\lim_{T\rightarrow\infty}\frac{1}{T}\int_0^T|\langle \xi|\psi(t)\rangle|^2 dt=
\sum_{j,k\,{\rm with }\, \theta_j=\theta_k }a_j a_k^* \langle\xi |\phi_j\rangle \langle\phi_k |\xi\rangle,
$$
where the asterisk $^*$ denotes the complex conjugate.
\end{theorem}
\begin{proof}
Since we have $U(t)=\sum_{j=1}^N\lambda_j(t)\, |\phi_j\rangle\langle\phi_j|$ and $|\psi(0)\rangle=\sum_{j=1}^N a_j\, |\phi_j\rangle$, the state of the system at time $t$ can be written as
$$
|\psi(t)\rangle=U(t)\,|\psi(0)\rangle=\sum_{j=1}^N\lambda_j(t)\, a_j\, |\phi_j\rangle
$$
using the orthonormality of the eigenstate basis.
Given a state $|\xi\rangle$, let us consider the projection of $|\xi\rangle$ on $|\psi(t)\rangle$. It holds
$$
\langle\xi |\psi(t)\rangle=\sum_{j=1}^N\lambda_j(t)\, a_j\,\langle\xi |\phi_j\rangle
$$
and therefore the probability of measuring the state $|\xi\rangle$ in $|\psi(t)\rangle$ is, for a given $t$,
\begin{eqnarray*}
&&|\langle\xi |\psi(t)\rangle|^2=\langle\xi |\psi(t)\rangle\,\langle\psi(t) |\xi\rangle=\\
&&=\sum_{j,k=1}^N \lambda_j(t)\lambda_k(t)^*\, a_ja_k^*\,\langle\xi |\phi_j\rangle\langle\phi_k |\xi\rangle.
\end{eqnarray*}
Now the time average can be written as
\begin{equation}
\frac{1}{T}\int_0^T|\langle\xi |\psi(t)\rangle|^2 dt=
\sum_{j,k=1}^N\frac{1}{T}\left(\int_0^T \lambda_j(t) \lambda_k(t)^* dt\right) a_ja_k^*\,\langle\xi |\phi_j\rangle\langle\phi_k |\xi\rangle.\label{eq:prob}
\end{equation}
Let us focus on the quantity 
$$\kappa=\lim_{T\rightarrow\infty}\frac{1}{T}\int_0^T \lambda_j(t) \lambda_k(t)^* dt=
\lim_{T\rightarrow\infty}\frac{1}{T}\int_0^T {\rm e}^{it(\theta_k-\theta_j)}dt.$$ 
If $\theta_j=\theta_k$, then the value of $\kappa$ is $1$. On the other hand, if $\theta_j\neq\theta_k$, the integral
$$
\int_{0}^T{\rm e}^{it(\theta_k-\theta_j)}dt=-\frac{i({\rm e}^{i(\theta_k-\theta_j) T}-1)}{\theta_k-\theta_j}
$$  
is uniformly bounded in modulus, and therefore $\kappa=0$. As a consequence, in the right-hand side of \eqref{eq:prob} the terms corresponding to pairs $(j,k)$ where $\theta_j\neq\theta_k$ are zero in the limit $T\rightarrow\infty$. From a physical viewpoint one might say that the interference between subspaces of different energy vanishes in the long-time average. The thesis follows immediately.
\end{proof}

Note that the limit (and hence the centrality scores) depend on the initial state and on the eigenstates of $H$, but not on the eigenvalues of $H$. 

\begin{corollary}\label{cor:cor1}
The limits in \eqref{eq:ourhub} and \eqref{eq:ourauthority} exist and are finite for every basis state $|j\rangle$. The associated hub and authority centrality measures are therefore well defined. 
\end{corollary}

\begin{corollary}\label{cor:cor2}
If the eigenvalues of $H$ are distinct, then \eqref{eq:ourhub} and \eqref{eq:ourauthority} simplify to
$$
C_{\rm hub}(j)=\sum_{i=1}^{2n} |a_i|^2 |\langle j|\phi_i\rangle|^2
$$
and
$$
C_{\rm auth}(j)=\sum_{i=1}^{2n} |a_i|^2 |\langle n+j|\phi_i\rangle|^2,
$$
for $j=1,\ldots,m$.
\end{corollary}

Let us now discuss the choice of the initial state $|\psi(0)\rangle$. This is a crucial part in the definition of our centrality measures, as highlighted by the results above. 

The most natural choice would seem to be the same as in \cite{Izaac}, i.e., a uniform occupation of the nodes:
\begin{equation}
|\psi(0)\rangle=\frac{1}{\sqrt{2n}}\sum_{k=1}^{2n} |k\rangle. \label{eq:initialuniform}
\end{equation}

The method that computes hub and authority scores as in \eqref{eq:ourhub} and \eqref{eq:ourauthority}, with the initial state defined as in \eqref{eq:initialuniform}, will be referred to as CQAu. (Here ``CQ'' stands for Continuous-time Quantum walk, ``A'' denotes the use of the adjacency matrix in the definition of the Hamiltonian, and ``u'' refers to the choice of a uniformly distributed initial state). 

For reasons that will be discussed below in this section and in Section \ref{sec:experiments}, we also consider the following choice for the initial state:

\begin{equation}
|\psi(0)\rangle=\frac{1}{\sqrt{\sum_{k=1}^{2n} d_k}}\sum_{k=1}^{2n} \sqrt{d_k}\,|k\rangle, \label{eq:initialweighted}
\end{equation}
where $d_k$ is the degree of the $k$-th node in the bipartite graph. Therefore, the initial occupation of the first $n$ nodes in the bipartite graph is proportional to the out-degrees of the nodes in the original directed graph, and the occupation of the last $n$ nodes is proportional to the in-degrees. The method that computes hub and authority  scores via \eqref{eq:ourhub} and \eqref{eq:ourauthority}, with the weighted initial state \eqref{eq:initialweighted}, will be referred to as CQAw. Despite the use of the same evolution operator, the two methods can lead to significantly different centrality rankings, as it will be seen in Section \ref{sec:experiments}.

Theorem \ref{th:Thm1} offers an interesting perspective on the choice of the initial state $|\psi(0)\rangle$ and on comparisons with other (classical) ranking methods where the eigenstructures of suitably chosen matrices play a relevant role. Let us assume for simplicity that in this discussion all matrices have simple, distinct eigenvalues. Consider first the definition of {\em subgraph centrality} (see e.g. \cite{EstradaHigham}) for an undirected $n$-node network: the score $c(i)$ of node $i$ is given by the $i$-th diagonal entry of the exponential of the (symmetric) adjacency matrix:
$$c(i)=[{\rm e}^A]_{i,i}=e_i^T {\rm e}^A e_i,$$
where $e_i$ is the $i$-th vector of the canonical basis in $\mathbb{R}^n$.
Let $A=V\Lambda V^T$ be the eigendecomposition of $A$, with $\Lambda={\rm diag} (\lambda_1,\ldots,\lambda_n)$. We assume here that the eigenvalues are in decreasing order. Then it holds ${\rm e}^A=\sum_{k=0}^n {\rm e}^{\lambda_k} v_k v_k^T$, where $v_1,\ldots, v_n$, the columns of $V$, are the eigenvectors of $A$. So the expression for $c(i)$ can be rewritten as
\begin{equation}
c(i)=\sum_{k=0}^n {\rm e}^{\lambda_k} (e_i^T v_k) (v_k^T e_i)=\sum_{k=0}^n {\rm e}^{\lambda_k} (V_{i,k})^2.\label{eq:expcentrality}
\end{equation}
In other words, we are looking at ${\rm e}^A$ as a linear combination of rank-one terms, each of them defined by one eigenvector and weighted with the exponential of the corresponding eigenvalue. The centrality $c(i)$ is the sum of the $i$-th diagonal entries of these rank-one terms and can be seen as a term ${\rm e}^{\lambda_1}(V_{i,k})^2$ associated with the dominant eigenvector, plus corrections with smaller coefficients depending on the other eigenvectors.

Now let us look at HITS. It has been already mentioned in Section \ref{sec:background} that this method computes the dominant eigenvectors of $AA^T$ and $A^TA$. Consider, as above, the eigendecomposition $AA^T=V\Lambda V^T$: then the hub score $h(i)$ of node $i$ can be written as
\begin{equation}
h(i)=\sum_{k=0}^n f(\lambda_k) (e_i^T v_k) (v_k^T e_i)=\sum_{k=0}^n f(\lambda_k) (V_{i,k})^2,\label{eq:HITShub}
\end{equation}  
where $f(x)$ is any function such that $f(\lambda_1)=1$ and $f(\lambda_j)=0$ for $j=2,\ldots,n$. A similar expression, obtained from the eigendecomposition of $A^TA$, holds for authority scores, and an analogous argument could be made for other ranking methods that ultimately rely on the computation of a dominant eigenvector. Generally speaking, we conclude that, for a large family of ranking methods, centrality scores are given by expressions like \eqref{eq:HITShub}, where $\{\lambda_j,v_j\}$ are the eigenpairs of the relevant matrix encoding the graph structure, and $f(x)$ is a suitably chosen function.

Now compare \eqref{eq:expcentrality} and \eqref{eq:HITShub} to Corollary \ref{cor:cor2}. In the CTQW setting, the coefficients $a_j$ of the initial state in the eigenvector basis clearly play a similar role to ${\rm e}^{\lambda_j}$ or $f(\lambda_j)$ for subgraph centrality and HITS, respectively. In other words, the choice of the initial state for a CTQW, where the relevant matrix is the evolution operator, corresponds to the choice of a weight function $f(x)$ in \eqref{eq:HITShub}.

\begin{remark}
It is well-known (see e.g., \cite{Ding}) that the HITS hub and authority scores correlate strongly with the out- and in-degrees of the nodes. This observation, combined with the discussion above, provides an additional motivation for the choice of the initial state $|\psi(0)\rangle$ as in \eqref{eq:initialweighted}.  
\end{remark}

\subsection{A method based on the Google matrix: CQG}

The construction of the Hamiltonian operator $H$ in \eqref{eq:rank1} and \eqref{eq:hamiltonian} is based on the adjacency matrix $A$. In principle, however, any matrix that encodes the structure of the graph can be employed to build a valid $H$ on the bipartite graph. For instance, in order to obtain rankings related to the classical PageRank, one can use the Google matrix $G$ from \eqref{eq:googlematrix} instead of $A$ and define
\begin{equation}
H_a=\left[\begin{array}{cc}
0 & G \\ G^T & 0
\end{array}\right].\label{eq:gmethod}
\end{equation}
In this setup, the matrix $G^T$ in the second block row, acting on the upper half of the state vector, describes a ``Google-like'' CTQW. Authority scores are extracted from the lower half of the limiting distribution. 

Note that the Google matrices used for classical and reverse PageRank are {\it not} the transpose of one another. Indeed, the correction of dangling nodes is not symmetrical, nor is the row normalization. In order to compute hub scores, we need to consider a different random walk described by a different Hamiltonian operator 
$$
H_h=\left[\begin{array}{cc}
0 & G_r \\ G_r^T & 0
\end{array}\right],
$$
where $G_r$ is the Google matrix obtained from the transpose of the adjacency matrix of the graph, as in reverse PageRank. As above, the hub scores are given by the lower half of the limiting distribution. Note that in this approach we need to run two separate CTQW to obtain hub and authority scores, whereas in CQAu and in CQAw a single random walk yields both.

The method described in this subsection, with the initial state defined as in \eqref{eq:initialuniform}, will be referred to as CQG in the following. 

Of course, the results of Theorem \ref{th:Thm1} and of Corollary \ref{cor:cor1} and \ref{cor:cor2} also hold in this case, since they apply to any unitary evolution.

\section{Numerical experiments}\label{sec:experiments}

We test the behavior of the three proposed CTQW-based centrality measures (CQAu, CQAw and CQG) on several examples of directed networks and compare them to other ranking algorithms: HITS, the method by Benzi et al. (BEK) \cite{BEK}, PageRank/reverse PageRank (PR) and, when computationally feasible, the discrete-time quantum PageRank presented in \cite{Paparo} (DQPR). In these experiments, HITS is always initialized with a uniform vector and implemented as a power method for the computation of the dominant eigenvalue of $AA^T$ and $A^TA$. DQPR is applied through an adaptation of Theorem 3.4 in \cite{Aharonov}, in combination with the speedup techniques outlined in \cite{Paparo}. 

In this paper we do not attempt a quantum implementation of the proposed centrality measures: all the numerical tests involve classical computations. The limits in \eqref{eq:ourhub} and \eqref{eq:ourauthority} are computed via the results of Theorem \ref{th:Thm1}. We also tried approximating these limits by computing integrals for increasing values of $T$ via numerical quadrature and stopping when the variation of the corresponding time average became smaller than a fixed tolerance. Note that, in the interest of speeding up these classical simulations, one may also apply fast methods for approximating the action of the matrix exponential on a vector or bilinear forms associated with a matrix function: see, e.g., \cite{AlMohy, Meurant}. The approach based on Theorem \ref{th:Thm1}, however, turned out to be faster.

All the methods were implemented in Octave and run on a laptop equipped with a 4-core Intel i7-7500U 2.70GHz processor, with 16 GB of RAM.

This section summarizes the main results obtained for the various methods. See Appendix \ref{app:A} for the complete data.

\subsection{Small graphs}

\begin{figure}
\centering
\includegraphics[width=0.6\textwidth]{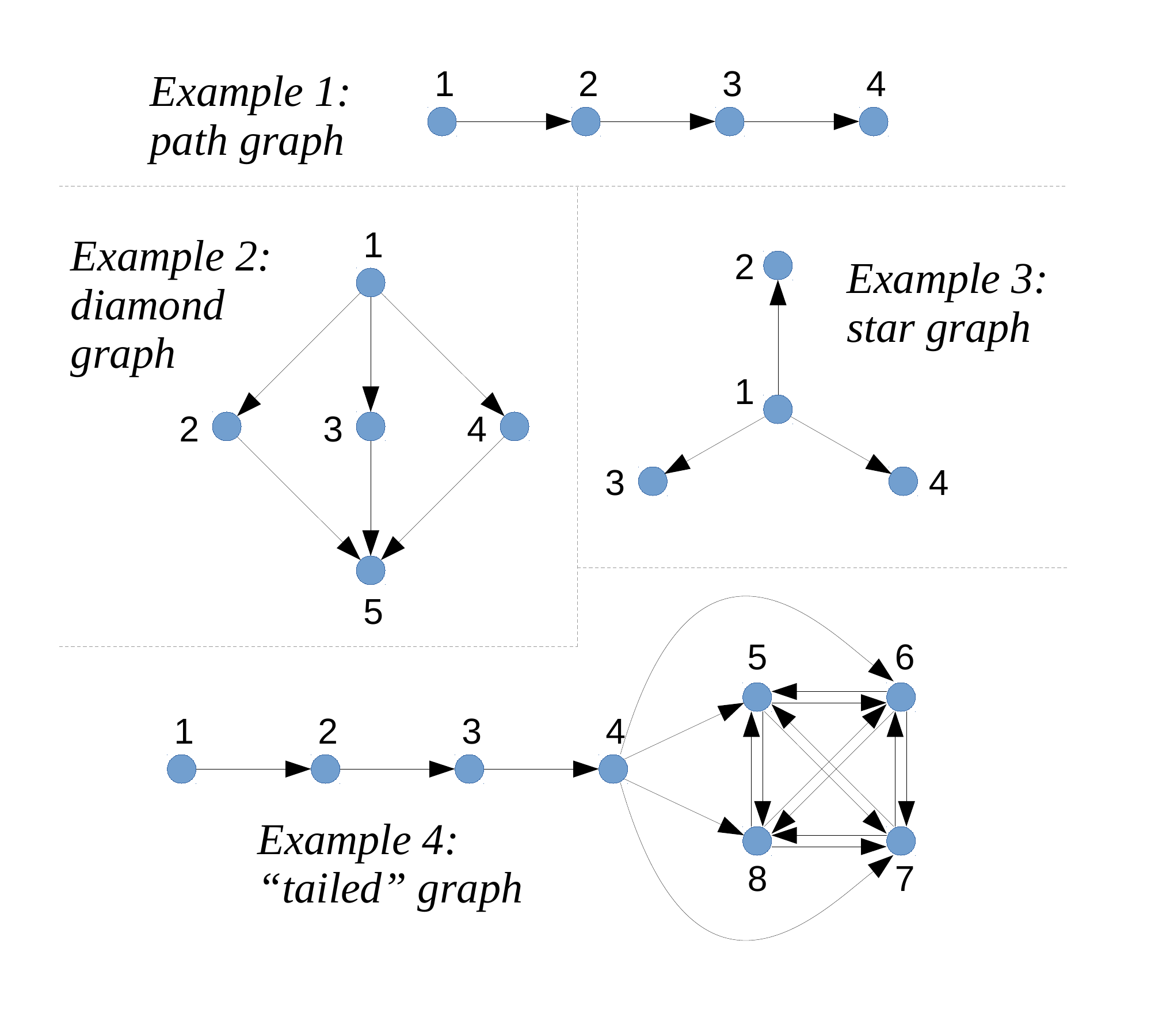}
\caption{Simple graphs used for testing the ranking methods (Examples 1 -- 4).}\label{fig:simplegraphs}
\end{figure}

\begin{figure}
\centering
\includegraphics[width=0.6\textwidth]{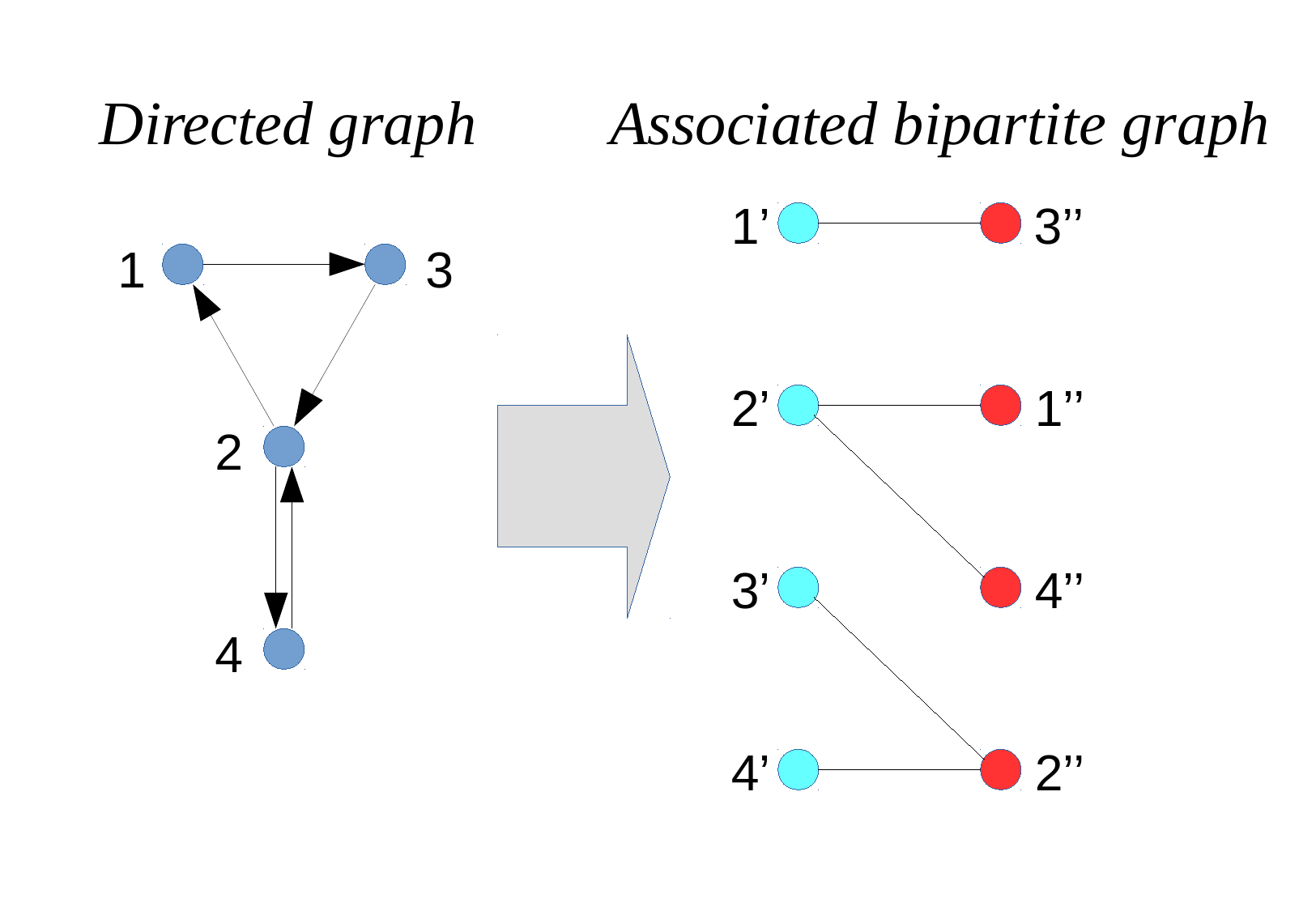}
\caption{Graph from Example 5. Despite the graph being strongly connected, the associated bipartite graph has 3 connected components.}\label{fig:BEK}
\end{figure}

In this subsection we compare ranking methods on a few small, simple graphs, where ``reasonable'' rankings are easily determined by inspection. These examples are useful to observe typical behaviours and pathologies of the tested methods. The examples are shown in Figures \ref{fig:simplegraphs} (examples 1-4) and \ref{fig:BEK} (example 5).
\begin{description}
\item{Example 1:} path graph, i.e., a simple chain of nodes.
\item{Example 2:} diamond graph, i.e., a graph with a ``main hub'' (node 1) with $n-2$ outgoing edges directed towards nodes $2,3,...,n-1$, and a ``main authority'' (node $n$) with $n-2$ incoming edges from nodes $2,3,...n-1$.
\item{Example 3:} star graph, i.e., a central hub (node 1)  connected to all the other nodes, and no other edge present.
\item{Example 4:} tailed graph, i.e., $n_1$ initial nodes connected as in a path graph, followed by $n_2$ nodes that form a complete subgraph. Node $n_1$ is connected to all the $n_2$ nodes in the complete subgraph.
\item{Example 5:} an example of strongly connected graph whose associate bipartite graph is not connected, taken from \cite{BEK}.
\end{description}

An interesting result emerging from these tests is the evidence of a pathology affecting almost all the quantum methods (including DQPR), except CQAw. When a graph contains ``no-hub'' nodes, i.e., nodes with zero out-degree, any sensible hub ranking should place them at the bottom. The same can be said for ``no-authority'' nodes (i.e., nodes with zero in-degree) in authority rankings. This expectation is satisfied 
 for all the classical methods and for CQAw, but not in the other quantum methods. CQAu fails in Example 2 (node $n$ comes second in the hub ranking; by symmetry, node $1$ comes second in the authority ranking) and in Example 3  (node 1 comes first as an authority). DQPR fails in Example 2, with the same behaviour as CQAu, and in Example 4, where node 1 is not last in the authority ranking. CQG fails on Example 4, where node 1 is slightly above nodes 2, 3 and 4 as an authority. This ``quantum pathology'' will have very noticeable consequences in an example of the next subsection, and affects especially CQAu.

In fact, it is the discovery of this anomalous behaviour in CQAu that led us to devise the method CQAw. In all the tests we performed, CQAw proved to be immune to the pathology, correctly ranking no-hub and no-authority nodes in the last positions\footnote{We also tried using a weighted initial vector for CQG, following the same idea as in CQAw. Such a modification does indeed improve the authority ranking in Example 4, but it has unpleasant consequences in simpler graphs -- such as the diamond and the star graphs -- where the original CQG method works correctly. So, the idea of a weighted initial vector was discarded; note also that CQG seems to be much less affected by the ``quantum pathology'' than CQAw.}.

The results of the numerical experiments also lead us to formulate the following observations:
\begin{itemize}
\item{} Recall that, in path graphs, HITS and PR display different behaviours. HITS singles out the last node (worse hub) or the first node (worse authority) and gives a tie for all the other nodes. PR, on the other hand, finds a strictly monotonic sequence of scores: increasing for authorities and decreasing for hubs. 

How do the quantum methods compare? It turns out that only DPQR behaves like PR, whereas all the other methods, including CQG, produce HITS-like rankings.   

\item{} Note that in our results HITS does not detect the top-ranked hubs and authorities in Examples 2 and 5, whereas the other methods do. HITS also fails to recognize the worse authority in Example 4, because node 1 has the same score as nodes 2, 3 and 4, although this is a less significant failure. However, note that Example 5 is a case where the hub and authority matrices $A^TA$ and $AA^T$ have a dominant eigenvalue of multiplicity 2; as a consequence, the computed eigenvectors -- and therefore the rankings -- depend on the initial vector and on implementation details. For a discussion of Example 5 see also \cite{BEK}.

\item{} Hub results for Example 4 highlight the difference between Google-like rankings, where tail nodes receive top scores, thus emphasizing the fact that they indirectly reach all the other nodes, and HITS-like rankings, which give more importance to the number of direct outgoing edges. Among the other methods, CQAw, BEK and -- curiously -- DQPR show a HITS-like behaviour, while CQG exhibits a PR-like behaviour.

\item{} DQPR gives an unexpected authority ranking in Example 4: the first two nodes in the ranking are 3 and 2, but they clearly do not have a relevant role as authorities in the graph structure. Moreover, node 1 is not ranked last as it should be.
\end{itemize}

In conclusion, from the tests on small graphs, CQAw appears to be immune to ``quantum pathologies'' and provides rankings that look overall quite similar to HITS. CQG is not immune to quantum pathologies, but does not seem to be strongly affected, either. Its rankings are mostly related to PR, with some exception, such as path graphs. CQAu, on the other hand, is too strongly affected by quantum pathologies to give reliable results on small graphs, but it usually gives the same result as CQAw at least on the top-ranked node.

\subsection{Artificial scale-free graphs}\label{sec:128nodes}
In this example, three 128-nodes graphs have been generated using the software NetworkX \cite{NetworkX}, according to the model developed in \cite{Bollobas} for the construction of scale-free graphs. The model requires three parameters: $\alpha$ and $\gamma$ determine the probability of adding new nodes according to the in-degree and out-degree distributions, respectively, whereas $\beta$ determines the probability of adding edges between pairs of existing nodes. Here the parameters have been chosen so as to generate three graphs with different properties. We examine the qualitative accordance of the top 10 nodes of the rankings, which are usually the most interesting, and the global correspondence of the rankings measured by Kendall's $\tau$ test.

\begin{figure}
\centering
\includegraphics[width=0.9\textwidth]{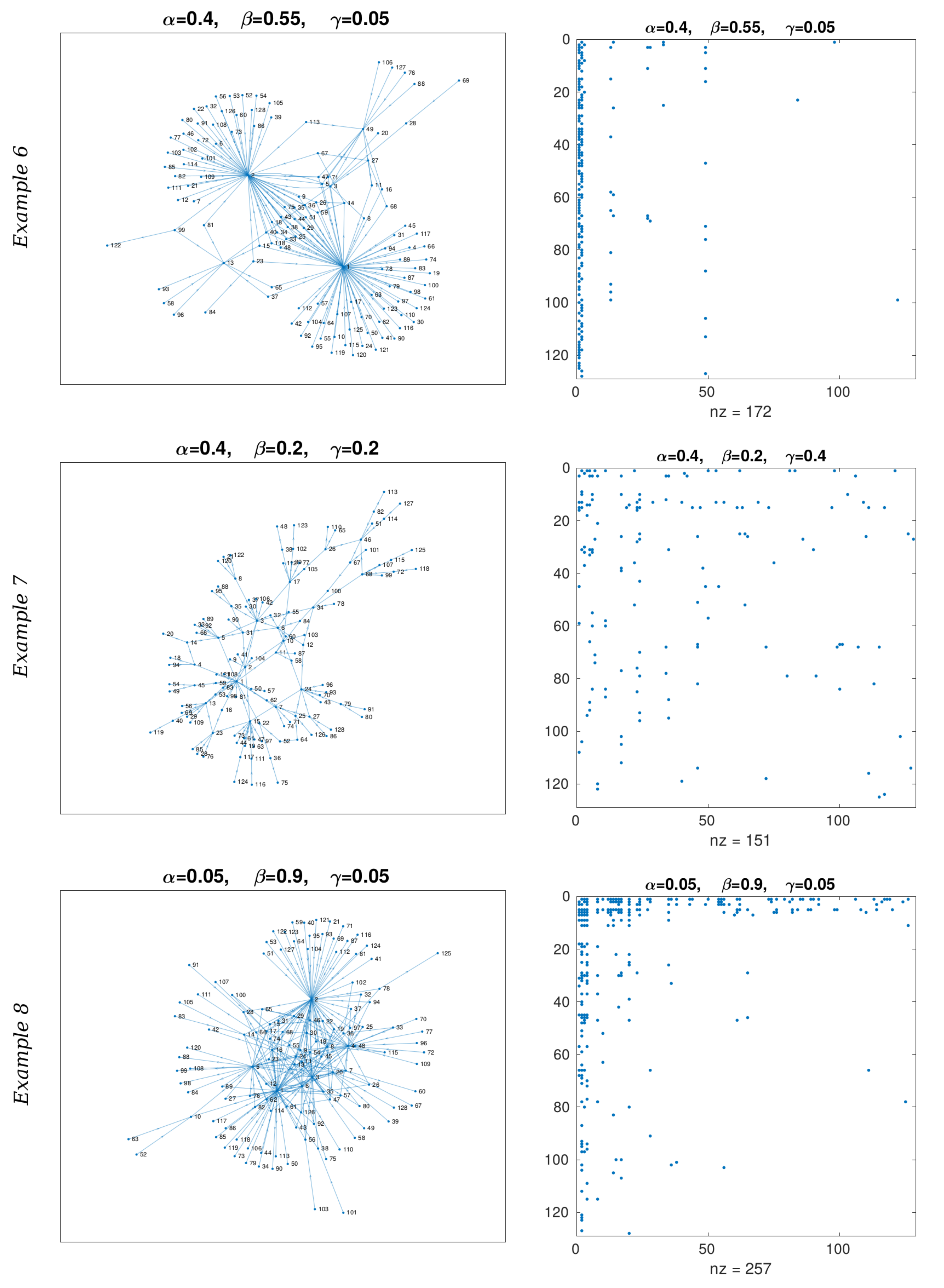}
\caption{Graphs used for testing the ranking methods (Examples 6, 7, 8). At the right of each graph the adjacency matrix is shown.}\label{fig:N128}
\end{figure}

The chosen examples are:
\begin{description}
\item{Example 6:} $\alpha = 0.4$, $\beta = 0.55$ and $\gamma = 0.05$. The graph shows a few strong authorities, among which nodes 1 and 2 are clearly dominant. There are also many (116) no-authority nodes. Hubs are present but weaker. The number of edges is 172.
\item{Example 7:} $\alpha = 0.4$, $\beta = 0.2$, $\gamma = 0.4$. The graph has a rather small number of edges, and no clearly dominant nodes. Hub and authorities are expected to be quite symmetrical. The number of edges is 151.
\item{Example 8:} $\alpha = 0.05$, $\beta = 0.9$, $\gamma = 0.05$. This graph has a larger number of edges (257), and a rather balanced hub/authority structure. Node 1 is likely to be the dominant hub, and node 2 the dominant authority.
\end{description}

The experiments show the effect of ``quantum pathologies'' on CQAu rankings: ``bubbles'' of ex-aequo no-hub or no-authority nodes tend to float up to high-ranked positions, disrupting the expected structure of the rankings and, in some cases, even reaching top positions. The most striking example is given by Example 6. In the hub ranking, the 9 no-hub nodes are ranked from 3rd to 11th, disrupting the top 10 ranking. (Incidentally, the same happens in DQPR, with the difference that the no-hub nodes are not ex-aequo).  In the authority ranking, the 116 no-authority nodes are placed 8th-121st, above 7 nodes that indeed {\it have} incoming edges. 

In the other two examples, the ``bubble'' does not reach the top 10 positions, but it is present anyway, as one can see by observing the overall ranking: in Example 7, positions 14th-74th for hubs and 19th-83rd for authorities; in Example 8, positions 14th-73rd for hubs and 16th-80th for authorities. So, when we look at the {\it global} rankings given by CQAu we cannot expect to find meaningful results. Kendall's $\tau$ between CQAu and any of the other methods, including CQAw, is typically quite low, sometimes even negative.

However, if one looks at the top-ranked nodes, avoiding the ``bubble'' in CQAu, there is usually an excellent agreement between CQAu and CQAw. So, the CQAu method may still have some interest if only a few top-ranked nodes are sought. Its advantage over CQAw lies in the preparation of the quantum initial state, which is easier for a uniform vector.

The CQAw method seems to provide fully satisfying results. It does not show any evident pathologies, and its rankings are generally in excellent agreement with HITS (and BEK) rankings. Kendall's $\tau$ between CQAw and HITS is always above 0.82. The top-ranked node and the set of top 3 nodes coincide in all the six computed rankings. The average superposition of the top 10 nodes (i.e., the number of nodes appearing in both top 10 lists given by CQAw and HITS, averaged over the six rankings) is 9.17.

CQG, on the whole, shows a good agreement with PR. Kendall's $\tau$ is almost always above 0.8 (except for hubs in Example 6, where it is 0.573), and the average superposition of the top 10 nodes is 7.58. However, some peculiarities can be observed. For instance, in the autority ranking for Example 6, node 84 is in second position, above node 2, whereas no other method places node 84 above the 10th position. In the hub ranking for Example 7, CQG is the only method that does not place node 15 first (it is 5th instead). In the authority ranking for Example 8 it fails to identify node 2 as the main authority, placing node 36 instead on top. We did not find a clear explanation for this behaviour.

\subsection{Larger real-life networks}\label{sec:realnetworks}

Here we test ranking algorithms on a few larger real-life graphs. As above, we examine the top 10 nodes by inspection and compare the global rankings via Kendall's $\tau$. In these examples we do not apply DQPR because of its high computational cost. 

The test networks are:
\begin{description}
\item{Example 9:} the metabolic network of the bacterium {\it E.~coli}, taken from companion data \cite{Barabasidata} to Barab\'asi's book \cite{Barabasibook}. Each node represents a metabolite, and, if there is a reaction where A is an input and B is a product, an edge A$\rightarrow$B is present. The network has 1039 nodes and 5802 edges.
\item{Example 10:} a network extracted from references between word categories in Roget's Thesaurus. The network is taken from the SuiteSparse Matrix Collection \cite{SparseTAMU}. It has 1022 nodes and 5075 links.
\item{Example 11:} a network of the web pages linking to www.epa.gov, which dates from 2006. The network is taken from the SuiteSparse Matrix Collection \cite{SparseTAMU}. The graph has 4472 nodes and 8965 links.
\end{description}

The results confirm the good agreement between CQAw and HITS: the top-ranked node is always the same, Kendall's $\tau$ is always above 0.66, and the average superposition of the top 10 nodes is 8.17. The agreement is especially good for Examples 9 and 10 ($\tau>0.82$ and never less than 9 coincident nodes), but a slight discordance is present in Example 11.

Note that in these examples the top 10 nodes for CQAu are in excellent agreement with the results of CQAw: on average, there are 9 coincident nodes. This result supports the claim that the method can be useful when looking for the top-ranked nodes in large graphs. However, Kendall's $\tau$ between the two methods is not always good, as it can fall below 0.5: this shows that the overall agreement is not as good as the agreement of the top 10 nodes.

The behaviour of CQG is sometimes puzzling. When looking at the first 10 nodes, we see a good agreement with the other methods in Example 9, a bad agreement in Example 10 (in both hub and authority rankings, only one of the top 10 nodes is present in PR rankings, and no nodes in HITS rankings) and, in Example 11, a very good agreement with PR for hubs (8 coincident nodes) and a very bad agreement for authorities (no coincident nodes). However, Kendall's $\tau$ paints a different picture, showing a poor agreement with PR in Example 9 ($\tau \approx 0.46$), a better agreement for Example 10 ($\tau \approx 0.57$, despite the results on the top 10 nodes) and a surprisingly good agreement ($\tau > 0.8$ for both hubs and authorities) in Example 11. No clear explanations were found for this apparently erratic behaviour.

\section{Conclusions}\label{sec:conclusions}

We have presented and tested three centrality methods for directed networks based on continuos-time quantum walks. The underlying common idea is the introduction of a Hermitian Hamiltonian matrix related to the undirected bipartite graph associated to the given network, therefore allowing for a unitary evolution in $2n$ dimensions. The results of the numerical tests lead to the following evaluations of the three proposed methods:

\begin{itemize}
\item{} The method CQAw gives completely satisfying results on all the tests. It seems to be the only quantum method among the tested ones (the three proposed, plus DQPR) that is immune to the typical ``quantum pathologies'' present in certain small graphs, a fact confirmed by further experiments performed on much larger diamond-, star- and tailed-graphs. The rankings for CQAw are strongly correlated with the rankings given by HITS. Moreover, the method provides hub and authority rankings in a single run. The only drawback we see is the potential difficulty in preparing an initial quantum state with the required mean occupation values, when implementing a quantum algorithm based on this method; note however that this drawback is common to other quantum ranking algorithms such as DQPR, and that efficient methods for the preparation of non-trivial quantum states are available (see e.g., \cite{GroverRudolph, HolmesMatsuura}).
\item{} The method CQAu, despite the use of the same evolution operator as CQAw, is the most affected by the quantum pathologies discussed above, and it fails on diamond- and star-graphs for any tested number of nodes. In graphs with a large number of no-hubs (or no-authority) nodes, these nodes are often ranked higher than weak hubs or weak authorities. Despite these problems, however, the method performed surprisingly well in identifying the top nodes in large graphs, in the tested examples: in the three tests with $n > 1000$ (Examples 9-11), the top 10 rankings turned out to be almost indistinguishable from CQAw. Also, for the three tests on 128-node graphs, the only unfavourable case -- when considering the top 10 nodes only -- is Example 6. Of course, quantum pathologies disrupt in any case the rankings below the top positions, and a good agreement of the whole ranking cannot be expected: this is confirmed by the usually poor results of Kendall's $\tau$ test when comparing CQAu with other methods. In conclusion, from these experiments it seems plausible that CQAu could be employed, for large networks and with the goal of identifying top-ranked nodes only, in the case where the preparation of a uniformly-occupied state proved to be advantageous enough w.r.t.~CQAw to justify a potentially less dependable performance.
\item The CQG method is the most puzzling of the three. It has the disadvantage of requiring two runs to obtain hub and authority rankings, as it needs two different evolution operators. It is only slightly affected by quantum pathologies. Usually it identifies top-ranked nodes that are compatible with other rankings (especially PR and DQPR), but its behaviour sometimes does not look consistent with the graph structure (see e.g. Example 8). On a same graph, it can give good results for hubs and bad results for authorities (Example 11). However, Kendall's $\tau$ between CQG and PR is usually high, showing that the two rankings, as a whole, are in fact closely related; but this property do not always translate to detection of the same top-ranked nodes. Such considerations discourage the use of CQG as presented here, although they suggest intriguing questions for future research.
\end{itemize}

\newpage

\section*{Appendix A}\label{app:A}

\begin{table}[!h]
\begin{center}
\caption{Hub scores for Example 1, with 4 nodes. Authority scores can be read by inverting the order of the nodes.}\label{table:expath}
\begin{tabular}{l l l l l l l l}
\hline
Node & CQAu & CQAw & CQG & HITS & BEK & PR & DQPR\\
\hline\\
1 & 0.13413 & 0.16505 & 0.15201 & 0.57735 & 1.54308 & 0.37015 & 0.38645\\
2 & 0.13413 & 0.16505 & 0.15201 & 0.57735 & 1.54308 & 0.29881 & 0.22468\\
3 & 0.13413 & 0.16505 & 0.15201 & 0.57735 & 1.54308 & 0.21489 & 0.22323\\
4 & 0.09760 & 0.00484 & 0.04396 & 0.00000 & 1.00000 & 0.11616 & 0.16564\\
\hline
\end{tabular}
\end{center}
\end{table}

\begin{table}[!h]
\begin{center}
\caption{Hub scores for Example 2, with 5 nodes. Authority scores can be read by exchanging node 1 and node 5.}\label{table:exdiamond}
\begin{tabular}{l l l l l l l l}
\hline
Node & CQAu & CQAw & CQG & HITS & BEK & PR & DQPR\\
\hline\\
1 & 0.20273 & 0.24431 & 0.26238 & 0.50000 & 2.91458 & 0.46835 & 0.45211\\
2 & 0.07000 & 0.08477 & 0.07029 & 0.50000 & 1.63819 & 0.14068 & 0.09117\\
3 & 0.07000 & 0.08477 & 0.07029 & 0.50000 & 1.63819 & 0.14068 & 0.09117\\
4 & 0.07000 & 0.08477 & 0.07029 & 0.50000 & 1.63819 & 0.14068 & 0.09117\\
5 & 0.08728 & 0.00139 & 0.02674 & 0.00000 & 1.00000 & 0.10962 & 0.27438\\
\hline
\end{tabular}
\end{center}
\end{table}

\begin{table}
\begin{center}
\caption{Hub scores for Example 3, with 4 nodes.}\label{table:exstarhub}
\begin{tabular}{l l l l l l l l}
\hline
Node & CQAu & CQAw & CQG & HITS & BEK & PR & DQPR\\
\hline\\
1 & 0.27227 & 0.49571 & 0.31268 & 1.00000 & 2.91458 & 0.54198 & 0.54484\\
2 & 0.07591 & 0.00143 & 0.06244 & 0.00000 & 1.00000 & 0.15267 & 0.15172\\
3 & 0.07591 & 0.00143 & 0.06244 & 0.00000 & 1.00000 & 0.15267 & 0.15172\\
4 & 0.07591 & 0.00143 & 0.06244 & 0.00000 & 1.00000 & 0.15267 & 0.15172\\
\hline
\end{tabular}
\end{center}
\end{table}

\begin{table}
\begin{center}
\caption{Authority scores for Example 3, with 4 nodes.}\label{table:exstarauth}
\begin{tabular}{l l l l l l l l}
\hline
Node & CQAu & CQAw & CQG & HITS & BEK & PR & DQPR\\
\hline\\
1 & 0.22752 & 0.00193 & 0.07733 & 0.00000 & 1.00000 & 0.20618 & 0.20676\\
2 & 0.09083 & 0.16602 & 0.14089 & 0.57735 & 1.63819 & 0.26461 & 0.26441\\
3 & 0.09083 & 0.16602 & 0.14089 & 0.57735 & 1.63819 & 0.26461 & 0.26441\\
4 & 0.09083 & 0.16602 & 0.14089 & 0.57735 & 1.63819 & 0.26461 & 0.26441\\
\hline
\end{tabular}
\end{center}
\end{table}

\begin{table}
\begin{center}
\caption{Hub rankings for Example 4, with $n_1=n_2=4$.}\label{table:extailhub}
\begin{tabular}{l l l l l l l l l}
\hline
CQAu & 4     & 1,2,3   &   &   & 5,6,7,8 &       &   &   \\
CQAw & 4     & 5,6,7,8 &   &   &         & 1,2,3 &   &   \\
CQG & 1,2,3 &         &   & 4 & 5,6,7,8 &       &   &   \\ 
HITS & 4     & 5,6,7,8 &   &   &         & 1,2,3 &   &   \\
BEK  & 4     & 5,6,7,8 &   &   &         & 1,2,3 &   &   \\
PR   & 1     & 2       & 3 & 4 & 5,6,7,8 &       &   &   \\
DQPR & 4     & 5,6,7,8 &   &   &         & 3     & 1 & 2 \\
\hline
\end{tabular}
\end{center}
\end{table}

\begin{table}
\begin{center}
\caption{Authority rankings for Example 4, with $n_1=n_2=4$.}\label{table:extailauth}
\begin{tabular}{l l l l l l l l}
\hline
CQAu & 5,6,7,8 &   &         & 2,3,4   &       &   & 1 \\
CQAw & 5,6,7,8 &   &         & 2,3,4   &       &   & 1 \\
CQG & 5,6,7,8 &   &         & 1       & 2,3,4 &   &   \\ 
HITS & 5,6,7,8 &   &         & 1,2,3,4 &       &   &   \\
BEK  & 5,6,7,8 &   &         & 2,3,4   &       &   & 1 \\
PR   & 5,6,7,8 &   &         & 4       & 3     & 2 & 1 \\
DQPR & 3       & 2 & 5,6,7,8 &         &       & 1 & 4 \\
\hline
\end{tabular}
\end{center}
\end{table}

\begin{table}
\begin{center}
\caption{Hub scores for Example 5. Authority scores can be read by exchanging node 1 and node 3.}\label{table:ex2BEK}
\begin{tabular}{l l l l l l l l}
\hline
Node & CQAu & CQAw & CQG & HITS & BEK & PR & QPR\\
\hline\\
1 & 0.07612 & 0.05714 & 0.12551 & 0.00001 & 1.54308 & 0.20916 & 0.13783 \\
2 & 0.20871 & 0.21788 & 0.25990 & 0.57735 & 2.17818 & 0.38694 & 0.41273 \\
3 & 0.10758 & 0.11249 & 0.05730 & 0.57735 & 1.58909 & 0.20195 & 0.25945 \\
4 & 0.10758 & 0.11249 & 0.05730 & 0.57735 & 1.58909 & 0.20195 & 0.18999 \\
\hline
\end{tabular}
\end{center}
\end{table}

\begin{table}
\begin{center}
\begin{threeparttable}
\caption{Hub rankings for Example 6. $G1$ and $G2$ indicates groups of ex-aequo nodes, listed in the notes.}\label{table:128authhub}
\begin{tabular}{l l l l l l l l l l l}
\hline
CQAu & 3 & 15      & {\it G1}\tnote{a} &    &       &       &    &    &       &               \\
CQAw & 3 & 5,47,71 &               &    & 26,59 &       & 15 & 23 & 25    & {\it G2}\tnote{b} \\
CQG & 1 & 99      & 23            & 3  & 2     & 69    & 25 & 67 & 26,59 &               \\ 
HITS & 3 & 5,47,71 &               &    & 15    & 26,59 &    & 25 & 23    & {\it G2}\tnote{b} \\
BEK  & 3 & 5,47,71 &               &    & 15    & 26,59 &    & 25 & 23    & {\it G2}\tnote{b} \\
PR   & 1 & 99      & 23            & 3  & 2     & 8     & 20 & 67 & 69    & 25            \\
DQPR & 3 & 98      & 28            & 33 & 122   & 84    & 27 & 14 & 13    & 49            \\
\hline
\end{tabular}
\begin{tablenotes}
\item[a] Nodes 13, 14, 27, 28, 33, 49, 84, 98, 122.
\item[b] Nodes 9, 18, 29, 34, 35, 36, 38, 40, 43, 44, 48, 51, 75, 118.
\end{tablenotes}
\end{threeparttable}
\end{center}
\end{table}

\begin{table}
\begin{center}
\caption{Authority rankings for Example 6.}\label{table:128authauth}
\begin{tabular}{l l l l l l l l l l l}
\hline
CQAu & 1 & 2  & 49 & 14  & 13 & 27 & 33 & \multicolumn{3}{l}{\it (116 ex-aequo)}   \\
CQAw & 1 & 2  & 49 & 14  & 13 & 27 & 33 & 3  & 98 & 84 \\
CQG & 1 & 84 & 2  & 122 & 49 & 13 & 27 & 3  & 28 & 14 \\ 
HITS & 1 & 2  & 49 & 13  & 14 & 27 & 33 & 3  & 28 & 84 \\
BEK  & 1 & 2  & 49 & 13  & 14 & 27 & 33 & 3  & 28 & 84 \\
PR   & 1 & 2  & 33 & 3   & 14 & 98 & 49 & 13 & 27 & 28 \\
DQPR & 1 & 2  & 3  & 49  & 33 & 13 & 14 & 98 & 27 & 28 \\
\hline
\end{tabular}
\end{center}
\end{table}

\begin{table}
\begin{center}
\caption{Kendall's $\tau$ for rankings for Example 6, for hubs (italics) and authorities (regular).}\label{table:128authtau}
\begin{tabular}{l l l l l l l l l l l}
\hline
& CQAu & CQAw & CQG & HITS & BEK & PR & DQPR\\
\hline
CQAu &       & {\it 0.378} & {\it 0.011} & {\it 0.342} & {\it 0.342} & {\it 0.427} & {\it 0.192} \\
CQAw & 0.200 &             & {\it 0.225} & {\it 0.935} & {\it 0.936} & {\it 0.645} & {\it -0.234}\\
CQG & 0.176 & 0.970       &             & {\it 0.166} & {\it 0.167} & {\it 0.573} & {\it 0.531} \\
HITS & 0.193 & 0.993       & 0.971       &             & {\it 0.999} & {\it 0.587} & {\it -0.299}\\
BEK  & 0.193 & 0.993       & 0.971       & 1.000       &             & {\it 0.588} & {\it -0.298}\\
PR   & 0.180 & 0.980       & 0.953       & 0.979       & 0.979       &             & {\it 0.110} \\
DQPR & 0.089 & 0.460       & 0.449       & 0.460       & 0.460       & 0.463       &             \\
\hline
\end{tabular}
\end{center}
\end{table}

\begin{table}
\begin{center}
\caption{Hub rankings for Example 7.}\label{table:128fewlinkshub}
\begin{tabular}{l l l l l l l l l l l}
\hline
CQAu & 15 & 1  & 3  & 13 & 31 & 10 & 68 & 12 & 26 & 51 \\
CQAw & 15 & 1  & 3  & 13 & 10 & 31 & 12 & 68 & 26 & 16 \\
CQG & 68 & 1  & 67 & 3  & 15 & 25 & 27 & 31 & 13 & 10 \\ 
HITS & 15 & 1  & 3  & 13 & 10 & 31 & 12 & 16 & 27 & 25 \\
BEK  & 15 & 1  & 3  & 13 & 10 & 31 & 68 & 12 & 27 & 26 \\
PR   & 15 & 1  & 13 & 3  & 68 & 45 & 27 & 26 & 2  & 25,79 \\
DQPR & 15 & 13 & 1  & 68 & 3  & 27 & 45 & 25 & 26 & 10 \\
\hline
\end{tabular}
\end{center}
\end{table}

\begin{table}
\begin{center}
\caption{Authority rankings for Example 7.}\label{table:128fewlinksauth}
\begin{tabular}{l l l l l l l l l l l}
\hline
CQAu & 2  & 6  & 1  & 24 & 5  & 46 & 17 & 90 & 23 & 3 \\
CQAw & 2  & 1  & 24 & 6  & 5  & 17 & 23 & 46 & 34 & 35 \\
CQG & 17 & 34 & 23 & 22 & 24 & 3  & 46 & 8  & 5  & 11 \\ 
HITS & 2  & 24 & 1  & 6  & 5  & 23 & 17 & 53 & 35 & 11 \\
BEK  & 24 & 2  & 6  & 5  & 17 & 1  & 23 & 46 & 11 & 35 \\
PR   & 17 & 24 & 5  & 3  & 23 & 6  & 1  & 8  & 11 & 2 \\
DQPR & 17 & 24 & 5  & 6  & 23 & 46 & 11 & 8  & 1  & 2 \\
\hline
\end{tabular}
\end{center}
\end{table}

\begin{table}
\begin{center}
\caption{Kendall's $\tau$ for rankings for Example 7, for hubs (italics) and authorities (regular).}\label{table:128fewlinkstau}
\begin{tabular}{l l l l l l l l}
\hline
& CQAu & CQAw & CQG & HITS & BEK & PR & DQPR\\
\hline
CQAu &       & {\it -0.157}& {\it -0.281}& {\it -0.315}& {\it -0.251}& {\it -0.205}& {\it -0.201}\\
CQAw & -0.084&             & {\it 0.749} & {\it 0.824} & {\it 0.888} & {\it 0.818} & {\it 0.718} \\
CQG & -0.138& 0.721       &             & {\it 0.653} & {\it 0.729} & {\it 0.871} & {\it 0.787} \\
HITS & -0.178& 0.897       & 0.693       &             & {\it 0.860} & {\it 0.684} & {\it 0.615} \\
BEK  & -0.120& 0.933       & 0.742       & 0.916       &             & {\it 0.751} & {\it 0.678} \\
PR   & -0.103& 0.790       & 0.856       & 0.735       & 0.789       &             & {\it 0.821} \\
DQPR & -0.109& 0.670       & 0.756       & 0.620       & 0.668       & 0.835       &             \\
\hline
\end{tabular}
\end{center}
\end{table}

\begin{table}
\begin{center}
\caption{Hub rankings for Example 8.}\label{table:128manylinkshub}
\begin{tabular}{l l l l l l l l l l l}
\hline
CQAu & 1 & 5  & 6 & 2  & 3 & 9 & 30 & 7  & 11 & 22 \\
CQAw & 1 & 5  & 6 & 2  & 3 & 9 & 7  & 11 & 30 & 45 \\
CQG & 1 & 66 & 5 & 78 & 2 & 6 & 11 & 29 & 7  & 3 \\ 
HITS & 1 & 5  & 6 & 2  & 3 & 9 & 7  & 11 & 30 & 45 \\
BEK  & 1 & 5  & 6 & 2  & 3 & 9 & 7  & 11 & 30 & 45 \\
PR   & 1 & 5  & 3 & 2  & 6 & 7 & 66 & 9  & 45 & 11 \\
DQPR & 5 & 1  & 6 & 3  & 2 & 7 & 66 & 11 & 9  & 29 \\
\hline
\end{tabular}
\end{center}
\end{table}

\begin{table}
\begin{center}
\caption{Authority rankings for Example 8.}\label{table:128manylinksauth}
\begin{tabular}{l l l l l l l l l l l}
\hline
CQAu & 2  & 4 & 3  & 20 & 1  & 24 & 16 & 13 & 8  & 17 \\
CQAw & 2  & 4 & 3  & 20 & 1  & 8  & 17 & 16 & 13 & 24 \\
CQG & 36 & 2 & 8  & 4  & 17 & 20 & 3  & 14 & 16 & 10 \\ 
HITS & 2  & 4 & 3  & 20 & 1  & 8  & 17 & 16 & 13 & 24 \\
BEK  & 2  & 4 & 3  & 20 & 1  & 8  & 17 & 16 & 13 & 24 \\
PR   & 2  & 4 & 20 & 1  & 3  & 14 & 17 & 8  & 16 & 28 \\
DQPR & 4  & 2 & 20 & 1  & 3  & 17 & 8  & 14 & 16 & 28 \\
\hline
\end{tabular}
\end{center}
\end{table}

\begin{table}
\begin{center}
\caption{Kendall's $\tau$ for rankings for Example 8, for hubs (italics) and authorities (regular face).}\label{table:128manylinkstau}
\begin{tabular}{l l l l l l l l}
\hline
& CQAu & CQAw & CQG & HITS & BEK & PR & DQPR\\
\hline
CQAu &       & {\it -0.193}& {\it -0.178}& {\it -0.209}& {\it -0.209}& {\it -0.248}& {\it -0.119}\\
CQAw & -0.192&             & {\it 0.795} & {\it 0.971} & {\it 0.971} & {\it 0.871} & {\it 0.708} \\
CQG & -0.182& 0.856       &             & {\it 0.787} & {\it 0.787} & {\it 0.813} & {\it 0.858} \\
HITS & -0.202& 0.989       & 0.853       &             & {\it 1.000} & {\it 0.858} & {\it 0.687} \\
BEK  & -0.202& 0.989       & 0.853       & 1.000       &             & {\it 0.858} & {\it 0.687} \\
PR   & -0.191& 0.877       & 0.872       & 0.875       & 0.875       &             & {\it 0.771} \\
DQPR & -0.135& 0.748       & 0.764       & 0.744       & 0.744       & 0.824       &             \\
\hline
\end{tabular}
\end{center}
\end{table}

\begin{table}
\begin{center}
\caption{Hub rankings for Example 9.}\label{table:metabolichub}
\begin{tabular}{l l l l l l l l l l l}
\hline
CQAu & 593 & 590 & 317 & 734 & 737 & 735 & 250 & 752 & 366 & 820 \\
CQAw & 593 & 590 & 317 & 734 & 737 & 250 & 735 & 366 & 820 & 563 \\
CQG & 593 & 590 & 317 & 734 & 737 & 276 & 219 & 736 & 382 & 991 \\
HITS & 593 & 317 & 590 & 734 & 737 & 250 & 366 & 563 & 735 & 820 \\
BEK  & 593 & 317 & 590 & 734 & 737 & 250 & 366 & 563 & 735 & 820 \\
PR   & 593 & 590 & 317 & 737 & 734 & 735 & 250 & 820 & 563 & 752 \\
\hline
\end{tabular}
\end{center}
\end{table}

\begin{table}
\begin{center}
\caption{Authority rankings for Example 9.}\label{table:metabolicauth}
\begin{tabular}{l l l l l l l l l l l}
\hline
CQAu & 590 & 250 & 820 & 831 & 593 & 279 & 365 & 735 & 739 & 734 \\
CQAw & 590 & 250 & 820 & 831 & 593 & 279 & 365 & 735 & 739 & 734 \\
CQG & 590 & 250 & 831 & 736 & 820 & 365 & 593 & 279 & 480 & 262 \\
HITS & 590 & 250 & 820 & 831 & 279 & 593 & 365 & 735 & 739 & 563 \\
BEK  & 590 & 250 & 820 & 831 & 279 & 593 & 365 & 735 & 739 & 563 \\
PR   & 590 & 820 & 250 & 593 & 736 & 831 & 279 & 317 & 365 & 735 \\
\hline
\end{tabular}
\end{center}
\end{table}

\begin{table}
\begin{center}
\caption{Kendall's $\tau$ for rankings for Example 9, for hubs (italics) and authorities (regular).}\label{table:metabolictau}
\begin{tabular}{l l l l l l l}
\hline
& CQAu & CQAw & CQG & HITS & BEK & PR\\
\hline
CQAu &       & {\it 0.631} & {\it 0.141} & {\it 0.542} & {\it 0.542} & {\it 0.200} \\
CQAw & 0.589 &             & {\it 0.291} & {\it 0.894} & {\it 0.894} & {\it 0.375} \\
CQG & 0.181 & 0.358       &             & {\it 0.253} & {\it 0.253} & {\it 0.469} \\
HITS & 0.446 & 0.834       & 0.278       &             & {\it 1.000} & {\it 0.336} \\
BEK  & 0.446 & 0.834       & 0.278       & 0.999       &             & {\it 0.336} \\
PR   & 0.300 & 0.541       & 0.461       & 0.468       & 0.468       &             \\
\hline
\end{tabular}
\end{center}
\end{table}

\begin{table}
\begin{center}
\caption{Hub rankings for Example 10.}\label{table:rogethub}
\begin{tabular}{l l l l l l l l l l l}
\hline
CQAu & 507 & 714 & 664 & 511 & 539 & 540 & 713 & 688 & 470 & 660 \\
CQAw & 507 & 714 & 664 & 511 & 539 & 540 & 713 & 470 & 688 & 660 \\
CQG & 629 & 945 & 392 & 110 & 103 & 186 & 9   & 213 & 374 & 44 \\
HITS & 507 & 714 & 664 & 511 & 539 & 540 & 713 & 470 & 660 & 469 \\
BEK  & 664 & 507 & 539 & 714 & 511 & 540 & 674 & 660 & 721 & 688 \\
PR   & 583 & 582 & 103 & 664 & 857 & 941 & 688 & 663 & 890 & 846 \\
\hline
\end{tabular}
\end{center}
\end{table}

\begin{table}
\begin{center}
\caption{Authority rankings for Example 10.}\label{table:rogetauth}
\begin{tabular}{l l l l l l l l l l l}
\hline
CQAu & 557 & 660 & 556 & 470  & 698  & 507 & 469 & 539 & 674 & 697 \\
CQAw & 557 & 660 & 556 & 470  & 698  & 507 & 469 & 539 & 674 & 697 \\
CQG & 93  & 651 & 566 & 675  & 171  & 856 & 220 & 914 & 668 & 267 \\
HITS & 557 & 660 & 470 & 556  & 698  & 507 & 469 & 674 & 539 & 486 \\
BEK  & 557 & 660 & 556 & 698  & 470  & 539 & 674 & 469 & 562 & 507 \\
PR   & 171 & 331 & 330 & 1001 & 1000 & 46  & 276 & 557 & 420 & 832 \\
\hline
\end{tabular}
\end{center}
\end{table}

\begin{table}
\begin{center}
\caption{Kendall's $\tau$ for rankings for Example 10, for hubs (italics) and authorities (regular face).}\label{table:rogettau}
\begin{tabular}{l l l l l l l}
\hline
& CQAu & CQAw & CQG & HITS & BEK & PR\\
\hline
CQAu &       & {\it 0.806} & {\it 0.266} & {\it 0.667} & {\it 0.739} & {\it 0.458} \\
CQAw & 0.819 &             & {\it 0.383} & {\it 0.834} & {\it 0.891} & {\it 0.509} \\
CQG & 0.239 & 0.317       &             & {\it 0.330} & {\it 0.429} & {\it 0.582} \\
HITS & 0.671 & 0.824       & 0.228       &             & {\it 0.818} & {\it 0.421} \\
BEK  & 0.746 & 0.889       & 0.332       & 0.809       &             & {\it 0.514} \\
PR   & 0.409 & 0.471       & 0.586       & 0.360       & 0.483       &             \\
\hline
\end{tabular}
\end{center}
\end{table}

\begin{table}
\begin{center}
\caption{Hub rankings for Example 11.}\label{table:EPAhub}
\begin{tabular}{l l l l l l l l l l l}
\hline
CQAu & 53 & 109      & 940,2796 &     & 7    & 61  & 75   & 12  & 120  & 92 \\
CQAw & 53 & 940,2796 &          & 61  & 7    & 75  & 77   & 120 & 3831 & 109 \\
CQG & 61 & 75       & 7        & 36  & 38   & 102 & 53   & 109 & 4571 & 395 \\
HITS & 53 & 940,2796 &          & 77  & 3831 & 61  & 686  & 120 & 12   & 38 \\
BEK  & 53 & 940,2796 &          & 75  & 61   & 77  & 3831 & 120 & 686  & 12 \\
PR   & 75 & 61       & 7        & 102 & 109  & 119 & 36   & 53  & 38   & 12 \\
\hline
\end{tabular}
\end{center}
\end{table}

\begin{table}
\begin{center}
\caption{Authority rankings for Example 11.}\label{table:EPAauth}
\begin{tabular}{l l l l l l l l l l l}
\hline
CQAu & 710  & 221  & 1321 & 3615 & 839  & 1152 & 2838 & 790  & 4700 & 599 \\
CQAw & 710  & 1321 & 221  & 839  & 2838 & 599  & 942  & 790  & 1262 & 956 \\
CQG & 2452 & 463  & 1030 & 1061 & 1881 & 1197 & 3222 & 2785 & 4730 & 3125 \\
HITS & 710  & 839  & 1321 & 956  & 1262 & 2233 & 2252 & 283  & 2262 & 2799 \\
BEK  & 710  & 1321 & 839  & 221  & 956  & 1262 & 2252 & 2233 & 2838 & 942 \\
PR   & 1247 & 2838 & 967  & 708  & 287  & 221  & 2175 & 1576 & 275  & 2799 \\
\hline
\end{tabular}
\end{center}
\end{table}

\begin{table}
\begin{center}
\caption{Kendall's $\tau$ for rankings for Example 11, for hub (italics) and authorities (normal).}\label{table:EPAtau}
\begin{tabular}{l l l l l l l}
\hline
& CQAu & CQAw & CQG & HITS & BEK & PR\\
\hline
CQAu &       & {\it 0.537} & {\it 0.672} & {\it 0.670} & {\it 0.671} & {\it 0.673} \\
CQAw & 0.467 &             & {\it 0.661} & {\it 0.667} & {\it 0.668} & {\it 0.665} \\
CQG & 0.500 & 0.696       &             & {\it 0.805} & {\it 0.805} & {\it 0.818} \\
HITS & 0.385 & 0.721       & 0.654       &             & {\it 0.826} & {\it 0.808} \\
BEK  & 0.383 & 0.726       & 0.635       & 0.796       &             & {\it 0.809} \\
PR   & 0.443 & 0.692       & 0.808       & 0.662       & 0.637       &             \\
\hline
\end{tabular}
\end{center}
\end{table}

\end{document}